\documentclass[12pt,reqno]{amsart}
\theoremstyle{plain}
\newtheorem{Th}{Theorem}[section]
\newtheorem{Lem}[Th]{Lemma}
\newtheorem{Prop}[Th]{Proposition}

\theoremstyle{definition}

\newenvironment{pf}{\noindent{\bf Proof.}}{\CQFD
}

\newcommand{\CQFD}
{%
\mbox{}%
\nolinebreak%
\hfill%
\rule{2mm}{2mm}%
\medbreak%
\par%
}

\renewcommand\r{\rho}

\newcommand{\R}{{\bf R}}

\newcommand{\C}{{\bf C}}

\def\Ho{\vbox{\offinterlineskip\hbox{\kern 3pt$\scriptstyle\circ$}
\kern
1pt\hbox{$H$}}}
\begin{document}
\title{ Monge-Amp\`ere Foliations for degenerate solutions}

\author{Morris KALKA and Giorgio PATRIZIO}
\date{}
\footnotetext{Much of this work was done while Kalka was
visiting the University in Florence and G. Patrizio  visited Tulane
University. The Authors thank the institutions for their support. G. Patrizio
acknowledges the support of the GNSAGA (INdAM).}

\maketitle
\begin{abstract}
We study the problem of the existence and the holomorphicity of the Monge-Amp\`ere foliation associated to a  plurisubharmonic solutions of the complex homogeneous Monge-Amp\`ere equation even at points of arbitrary degeneracy. We obtain good results  for real analytic unbounded solutions. As a consequence we also provide a positive answer to a question of Burns  on homogeneous polynomials whose logarithm satisfies the complex Monge-Amp\`ere equation and  we obtain a generalization the work of P.M. Wong  on the classification of complete weighted circular domains.\\
\vskip0,1cm
\noindent
\vskip0,1cm 
\noindent
\end{abstract}

\section {Introduction}
\setcounter{section}{1}
\setcounter{equation}{0}

Let $M$ be a Stein manifold of dimension $n>1$ and suppose that there is  a continuous plurisubharmonic exhaustion
 $\rho \colon M \to [0,+\infty)$ 
 of class $C^\infty$ on 
$M_{*}=\{\rho >0\}$
with $d\rho \not= 0$ on $\{\rho >0\}$ such that 
 on $M_{*}$  the function $u=\log \rho$ satisfies the complex
homogeneous Monge-Amp\`ere equation,
\begin{align} \label{eqMA} (dd^c u)^n = 0.
\end{align}
\\
On the open set $P$ of $M_{*}$  where $dd^c\rho >0$, it is well known that 
\({\rm Ann}(dd^cu)\) is an integrable distribution of complex rank $1$ which is generated  by the complex gradient  $Z$ of
$\rho$ (with respect to the K\"ahler metric with potential \(\rho\))which, in local coordinates, is given by
\begin{align}
Z=\sum_{\mu,\nu}\rho^{\mu\bar{\nu}}\rho_{\bar\nu}\frac{\partial}{\partial
z^{\mu}}.
 \end{align}
Here and throughtout the paper lower indices denote derivatives and 
$(\rho^{\mu\bar{\nu}})$ is the matrix  defined by 
the relation $\sum_{\nu=1}^{n} \rho^{\mu\bar{\nu}}\rho_{\alpha\bar{\nu}}=
\delta_{\alpha \mu}$.
\\

It is natural to ask whether the vector field $Z$ and therefore the foliation by complex curves defined by the distribution \({\rm Ann}(dd^cu)\) -- known as the {\sl Monge-Amp\`ere foliation associated to}
$u$ -- extends throughout the degeneracy set $M_{*}\setminus P$.
This is the case in many natural examples. Let 
$\rho\colon {\bf C}^n\to [0,+\infty)$ be a continuous plurisuharmonic exhaustion, smooth on ${\bf C}^n_{*}=\{\rho >0\}$ satisfying
for suitable positive real numbers $c_{1},\dots,c_{n}$, the following ``homogeneity'' condition for $z=(z^{1},\dots,z^{n})\in{\bf C}^n$:
\begin{align} \label{genminkowski}
\rho(e^{c_{1}\lambda}z^{1},\dots,e^{c_{n}\lambda}z^{n})=
\vert e^{\lambda}\vert^{2}\rho(z).
\end{align}
The sublevel set of such an exhaustion have been referred in \cite{Wong} as  {\sl generalized weighted circular domains} a class of domains which   includes the complete circular domains. 
Using \eqref{genminkowski} it is easy to show that  $u=\log \rho$ satisfies the complex homogeneous Monge-Amp\`ere equation \eqref{eqMA}. In this case the complex gradient vector field extends on all 
$ {\bf C}^n$ and in fact it may  be computed for $z=(z^{1},\dots,z^{n})\in{\bf C}^n$: 
\begin{align}
Z=\sum_{\alpha} c_{\alpha}z^{\alpha}\frac{\partial\;}{\partial z^{\alpha}}.
 \end{align}
so that the leaves (i.e. the integral submanifolds of $Z$) are exactly the orbits of the action 
$z=(z^{1},\dots,z^{n})\in{\bf C}^n 
\mapsto (e^{c_{1}\lambda}z^{1},\dots,e^{c_{n}\lambda}z^{n})$
(see again (\cite{Wong}) for details).

In this paper we look for condition which grants the extendability of gradient vector field $Z$ and consequentely of the Monge-Amp\`ere foliation.
Indeed, if the exhaustion $\rho$ is real analytic on $M_{*}$
 we are able to extend  the bundle \({\rm Ann}(dd^cu)\)  as an integrable complex line bundle on
\(M_{*}\) generated by an extension of the complex gradient $Z$. Furthermore we prove that the extended vector field $Z$, and therefore  the associated Monge-Amp\`ere foliation, is holomorphic.
The same result is shown to be true if the exhaustion $\rho$ is of class $C^{\infty}$ on $M_{*}$ but the degeneracy is not too high.
As a consequence, a repetition of the arguments in our previous work 
(\cite{Kalka-Patrizio1},\cite{Kalka-Patrizio2}) 
 allows us to answer in full generality the  question of Burns (\cite{Burns}) regarding the characterizations of plurisuharmonic homogeneous polynomial and  to extend the results of P.M. Wong (\cite{Wong}) about classifying
complete weighted circular domains   assuming 
real analyticity of \(\rho\) and connectedness  of the set \(M_*\).
 In complex dimension \(n=2\), in  \cite{Kalka-Patrizio2} we were able to prove the same result if \(\rho\) is of class \(C^{\infty}\)on \(M_{*}\) and satisfies a finite type condition. For \(n>2\) it turns out that real analyticity is a crucial element of our arguments.

\section {Extension and holomorphicity of the foliation on low-degeneracy set}
\setcounter{equation}{0}

Throughout this section we shall make the following assumptions:
\\ \\
$(A1)$ $M$ is a Stein manifold of dimension $n>1$;
\\ \\
$(A2)$ $\rho \colon M \to \to [0,+\infty)$ is
continuous plurisubharmonic exhaustion of class $C^\omega$ 
with $d\rho \not= 0$ on $M_{*}=\{\rho >0\}$;
\\ \\
$(A3)$  on $M_*$  the function $u=\log \rho$ satisfies the complex homogeneous Monge-Amp\`ere equation,
\begin{align}(dd^c u)^n = 0.
\end{align}
\\

We start recalling the  notion of hypersurface of finite type in the
sense of Kohn \cite{Kohn} and Bloom-Graham \cite{BloomGraham}.
Let $S$ be a (real) hypersurface in a complex manifold $M$ and $p\in S$. Let $\rho
\colon U \to \R$ be a local defining function of $S$ at $p$, i.e.
$p\in U\cap S =\{z\in U \mid \rho(z)=0\}$ and $d\rho\neq 0$ on $U\cap S$.
Let  ${\mathcal L}^k$ denote the module over $C^{\infty}(U)$ functions generated by
all vector fields $W$ on $U$ in the holomorphic tangent bundle of $S$ (i.e.
with $\partial \rho(W)=0$ on $S$), their conjugates and their  brackets of length
at most $k$.  Then ${\mathcal L}^0$ is the module of vector fields on $U$ spanned
by the vector fields on $U$ holomorphically tangent to $S$ and their conjugates
and, for $k>0$,  ${\mathcal L}^k$ is the
the module of vector fields on $U$ spanned by brackets $[V,W]$ with
$V\in {\mathcal L}^{k-1}$ and $W\in{\mathcal L}^0$. The point $p\in S$ is  of {\sl type
m}  if $\partial \rho(W_{|p})=0$ for any $W\in {\mathcal L}^k$, for all $k=0,\dots, m-1$, and 
and there exists $Y\in {\mathcal L}^m$ with $\partial \rho(Y_{|p})\neq 0$. We say
that $S$ is of {\sl finite type} if for every $p\in S$ it is of type $m$ at
$p$, where $m$ may depend on $p$.
\\ \\

If the exhaustion $\rho$ satisfies (A2), it is  well known that for any $c>0$, the level hypersurface  $\{\rho=c\}$ is compact, real analytic, pseudoconvex and  of finite type. In fact, $d\rho\neq0$, $\rho$ is plurisubharmonic and if, for some $c_{0}>0$, $p$ were a point of  the real analytic hypersurface  $S=\{\rho =c_{0}\}$ which is not of finite type, there  would   be a complex  variety  of positive dimension through $p$ lying on S. This cannot happen  because is $S$ compact (see   \cite{DiederichFornaess}).
\\ \\

We shall use the following notations:

$$P=P_{n}=\{p\in M_{*}\mid (dd^{c}\rho)^{n}\neq0\},$$
$$P_{n-1}=\{p\in M_{*}\mid (dd^{c}\rho)^{n-1}\neq0\}.$$

We have the following:

\begin{Lem} Suppose that the assumptions $(A1)$, $(A2)$ and $(A3)$ hold. Then $P\subset P_{n-1}$ and both are open dense sets in 
$M_{*}$ (and hence also in $M$!).
\end{Lem}
\begin{pf}  Both $P$ and $P_{n-1}$ are open since they have closed  complements. If  $P_{n-1}$  were not a  dense set, there would be an open set $A\subset M_{*}$ in its complement. But then for any $p\in A$ the Levi form at $p$ restricted to the holomophic tangent space to the level set $S$ of $\rho$ through $p$ would have at least one zero eigenvalue. Hence  $S\cap A$ would  
 contain  complex varieties of positive dimension which, as we observed above, cannot happen. On the other hand, since
\begin{align} \label{ddcu}
 dd^cu=\frac{dd^c \r}{\r}-\frac{d\r\wedge d^c\r}{\r^2},
\end{align}
then from (A3) it follows that on $M_{*}$
\begin{align}
 \rho (dd^c\rho)^{n}=n (dd^c\rho)^{n-1}\wedge d\rho\wedge d^{c}\rho.
\end{align}
Suppose there were an open set $B\subset M_{*}$
on which $(dd^c\rho)^{n}=0$. Then let $p\in B$ and let $S$ be the level set of $\rho$ through $p$. If $V,Z\in T_{p}^{1,0}(M_{*})$ with
$d\rho_{p}(V)=0$ and $d\rho_{p}(Z)\neq0$, then
\begin{align}
 0= \rho (dd^c\rho)^{n}(V,\bar V,\dots,V,\bar V,Z,\bar Z) =
 n (dd^c\rho)^{n-1}(V,\bar V,\dots,V,\bar V)\vert d\rho (Z)\vert^{2}.
\end{align}
This allows us to argue as above and conclude that $S\cap B$ contains complex varieties,  which
 cannot happen.  
 \end{pf}
 
\vskip0.3cm
On $P$, $dd^c\rho >0$ and under the assumption (A3) it follows that 
 \begin{align}
{\rm Ann}(dd^cu) =\{W\in T^{1,0}(M_{*})\mid
				dd^cu(W,\bar V)=0 \hskip.3cm \forall v\in T^{1,0}(M_{*})\}
			 \end{align}
is an integrable distribution of complex rank $1$ which is generated  by the complex gradient  $Z$ of
$\rho$ (with respect to the K\"ahler metric with potential \(\rho\))
(see for instance \cite{Stoll} or \cite{Burns}  or  \cite{Wong} for details). A computation  in local coordinates shows that $Z$ is given by
\begin{align}\label{zeta}
Z=\sum_{\mu,\nu}\rho^{\mu\bar{\nu}}\rho_{\bar\nu}\frac{\partial}{\partial
z^{\mu}}.
 \end{align}
 
  It is well known (again see   \cite{Stoll} or \cite{Burns}  or  \cite{Wong} for details) that, wherever it is defined, the vector field $Z$ is a non vanishing section of the line bundle \({\rm Ann}(dd^cu)\) and the Monge-Amp\`ere equation \eqref{eqMA} is
equivalent to 
\begin{align} \label{Zrho}
Z(\rho)=\rho.
\end{align}

In this section we aim to extend the vector field 
$Z$ and the distribution ${\rm Ann}(dd^cu)$ to $P_{n-1}$. Consequently the open set $P_{n-1}$ will be foliated by complex curves, the flow of the extended vector  field $Z$. Using  the unboundness of $\rho$ in a crucial way, we  shall also show that the corresponding foliation of $P_{n-1}$ thus obtained is holomorphic.

To prove that \({\rm Ann}(dd^cu)\) extends as a complex line bundle on 
$P_{n-1}$ and it is an integrable distribution spanned by a suitable
extension of $Z$,   we first need  the following remark:

\begin{Lem} Suppose that  $(A1)$, $(A2)$ and $(A3)$ hold. For any $p\in P_{n-1}$ there exist a coordinate neighborhood  $U\subset \{\rho >0\}$  such that:
\\ $(i)$  there are real functions 
$\lambda_{1},\dots,\lambda_{n-1}$ on $U$ such that at every  point
of $q\in U$ the values $\lambda_{1}(q),\dots,\lambda_{n-1}(q)$ 
are the eigenvalues of the restriction of the Levi form  of $\log \rho$ to ${\rm Ker} \partial\rho(q)$. The coordinates on $U$ may be chosen so that if $\lambda_{j}(q),\lambda_{k}(q)\neq0$ for $q\in U$  and $j\neq  k$ then $\lambda_{j}(q)\neq\lambda_{k}(q)$ and the functions $\lambda_{1},\dots,\lambda_{n-1}$ are  of class $C^{\omega}$;
\\ 
$(ii)$  there are linearly independent vector fields
$L_{1},\dots,L_{n-1}\in {\rm Ker} \partial\rho$ on $U$ which are linear analytic and form a diagonalizing basis for the hermitian form defined by the restriction of the Levi form of $\log \rho$ to 
$ {\rm Ker} \partial\rho(q)$ at every point
 $q\in U$ so that for each $j$ the vector $L_{j}(q)$ is the eigenvector corresponding to $\lambda_{j}(q)$. 

\end{Lem} \label{diagogonal}
\begin{pf} At $p\in P_{n-1}$, and hence for every $q$ in a neighborhood $U$ of $p$, with respect to any coordinate system,
 at most one eigenvalue of $\rho_{\mu\bar \nu}$ is zero. 
 Therefore the restriction of the Levi form of $\rho$ to the holomorphic tangent space ${\rm Ker} \partial\rho(q)$ to the level set $S_{q}=\{\rho=\rho(q)\}\cap U$ through $q$ has at most one zero eigenvalue
 The same happens for the restriction of the Levi form of any other defining function of $S_{q}$ at $q$, thus, in particular for the restriction of the Levi form of  $\log \rho$.
 Thus, by suitably rescaling coordinates, it may be assumed that the eigenvalues of the restriction of the Levi form of $\log\rho$ are pairwise distinct and therefore $(i)$ follows. If $(i)$ holds, then (ii) is straightforward. 
\end{pf}


\bigskip
\bigskip\bigskip
In a coordinate  neighborhood $U\subset \{\rho >0\}$ as in  Lemma $2.1$,
the module ${\mathcal L}^k$  
 is generated over  $C^{\infty}$ functions by $L_{1},\dots,L_{n-1}$ and $\bar L_{1},\dots, \bar L_{n-1}$ and their  brackets of length at most
$k$.  Since the level sets of $\rho$ are hypersurfaces of finite
type, if $p\in U$, then there exists $k$ such that ${\mathcal L}^k_{|p}$ is the full complex tangent space $T^{\C}_{p}(M)$.


Recall that on  $P$ we have that $\partial\rho (Z)=Z(\rho)=\rho$ and hence  the complex gradient of $\rho$ is
transverse to holomorphic tangent bundle of any level set of $\rho$.
We shall extend  the bundle ${\rm Ann}(dd^cu)$ to all $P_{n-1}$ by suitably choosing  the ``missing direction'' recovered by the finite type property at the weakly pseudoconvex   points in $P_{n-1}\setminus P$.

We can  prove

\begin{Prop} Suppose that the assumptions $(A1)$, $(A2)$ and $(A3)$ hold. The complex gradient
$Z$ defined on $P$  extends to a  non--zero $C^{\infty}$ vector field
on   $P_{n-1}$. The extension of the vector field $Z$  generates an integrable complex line bundle ${\mathcal A}$. In particular, the
Monge-Amp\`ere foliation defined by  ${\rm Ann}(dd^c u)$ on $P$ extends to a foliation defined by the distribution ${\mathcal A}$ on   $P_{n-1}$.
 \end{Prop}
\begin{pf} Let $p\in P_{n-1}\setminus P$ and  $\lambda_{1},\dots,\lambda_{n-1}$ be the functions which provide the eigenvalues of the Levi form of $\log\rho$ on ${\rm Ker} \partial \rho$. Let also  $L_{1},\dots,L_{n-1}\in {\rm Ker} \partial\rho$ be the vector fields of class $C^{\omega}$ defined in  Lemma  $2.\ref{diagogonal}$ on an open coordinate set
$U\subset\{\rho >0\}$ containing $p$. 
\par
Since  the boundary of each sublevel set
of $\rho$ is of finite type,  for some positive integer $m$ there exists a
$C^{\infty}$  nonzero vector field $Y\in {\mathcal L}^{m}$ with $\partial \rho(Y)\neq
0$ at $p$ and hence in a neighborhood of $p$ (which we may assume is $U$).
Since $(Z-\bar Z)_{|q}$, ${L_{1}}_{|q},\dots,{L_{n-1}}_{|q}$ and 
$ {\bar {L_{1}}}_{|q},\dots, {\bar {L_{n-1}}}_{|q}$ span the tangent space
to the level set of $\rho$ through any $q\in U\cap P$,
 it follows that
\begin{align}\label{defY}
Y=\phi(Z-\bar Z)+ A_{1}L_{1}+\dots+ A_{n-1}L_{n-1} +
B_{1}\bar L_{1}+\dots+B_{n-1}\bar L_{n-1}
\end{align}
for suitable functions $\phi, A_{1},\dots,A_{n-1},
B_{1},\dots,B_{n-1}$ of class $C^{\infty}$ on $U\cap P$.
If we define
\begin{align}\label{defV}
V=\frac{1}{2}(Y-iJY),
\end{align}
then on $U\cap P$ we have
\begin{align}\label{computeV}
V=\frac{1}{2}(Y-iJY)=\phi Z + A_{1}L_{1}+\dots+ A_{n-1}L_{n-1}.
\end{align}
We  now  compute and study the functions $\phi$ and
$A_{1},\dots,A_{n-1}$. On  $U$ the function $\partial \rho(V)$ is smooth and non zero.
On the other hand on $U\cap P$, using (\ref{Zrho}) and (\ref{computeV})
\begin{align}
\partial \rho(V)=\phi\partial \rho(Z)=\phi\rho
\end{align}
so that
\begin{align}
\phi=\frac{\partial \rho(V)}{\rho}
\end{align}
extends as a non zero function of class $C^{\infty}$ on all $U$.
\\
On $U\cap P$, consider the form
\begin{align}\label{defineOmega}
\Omega=\frac{ (dd^{c}u)^{n-1}\wedge d\rho\wedge d^{c}\rho}
{\lambda_{1}\dots\lambda_{n-1}}.
\end{align}
The form $\Omega$ is a volume form on  $U$   which allows to compute and study the functions $A_{1},\dots,A_{n-1}$. First of all notice that
if $W$ is any vector field on $U$ such that $L_{1},\dots,L_{n-1},W$
span $T^{1,0}_{q}(M)$ at every point $q\in U$ then necessarily
$d\rho(W)\neq0$. If $\theta_{1},\dots,\theta_{n}$ is the coframe dual 
to $L_{1},\dots,L_{n-1},W$ on $U$, then
\begin{align*}
\Omega
(L_{1},\bar L_{1},\dots,L_{n-1} \bar L_{n-1}, W,\bar W) &= 
\frac{i}{2\pi} 
\frac{\prod^{n-1}_{j=1}dd^{c}u(L_{j},\bar L_{j})}
{\lambda_{1}\dots\lambda_{n-1}}
\vert d\rho(W)\vert^{2} \\&= 
\frac{i}{2\pi}
\vert d\rho(W)\vert^{2}.
\end{align*}
As a consequence
\begin{align}
\Omega=\frac{i}{2\pi}
\vert d\rho(W)\vert^{2}\theta_{1}\wedge\bar\theta_{1}
\wedge\dots\wedge\theta_{n}\wedge\bar\theta_{n}
\end{align}
is a $C^{\infty}$  non--zero volume form on $U$.
On the other hand we can compute for each $j=1,\dots,n-1$:
\begin{align}
\Omega
(L_{1},\bar L_{1},\dots,
V,\bar L_{j}\dots,
L_{n-1} \bar L_{n-1}, W,\bar W)=\frac{i}{2\pi}
\vert d\rho(W)\vert^{2}A_{j}
\end{align}
from which it follows that $A_{j}$ extends 
as   a $C^{\infty}$ functions  throughout $U$. 
\par
Therefore the complex gradient
$Z$ extends as a  $C^{\infty}$  vector field everywhere on $U$
by setting $Z=\frac{1}{\phi}(V-A_{1}L_{1}-\dots -A_{n-1}L_{n-1})$.
The rest of the statement is now obvious.
\end{pf}

\noindent
{\bf Definition.}  We call the foliation defined by  ${\mathcal A}$ the {\sl extended Monge-Amp\`ere foliation}. 
\\
\\

Before moving on, we need more precise information about
the extended  Monge-Amp\`ere foliation. In particular we need to show that a leaf starting at a point in the open set $P$ lays entirely in $P$.
In this step we use the use the hypothesis real analyticity.
\\

We start recalling a few known facts proved by Bedford an Burns (\cite{Bedford-Burns}). 
 Let $p\in P$ and let $\ell_{p}$  be the leaf through $p$. According to \cite{Bedford-Burns} the leaf  $\ell_{p}$ extends also through the set where $dd^{c}u$ has rank strictly less than $n-1$ to yield a ``complete'' leaf that we still call $\ell_{p}$. 
The vector field $Z$, originally defined on the  points of $\ell_{p} \cap P$, extends on all $\ell_{p}$ and, by continuity satisfies $d\rho(Z)=\rho$. Hence in particular it is non zero. As a consequence  the vector field $Z$ is defined on the union of all the leaves 
passing through points of $P$. Furthermore we recall that in \cite{Bedford-Burns} it is shown that if $p\in P$ and $\ell_{p}$ is the extended leaf through $p$ thus obtained, the set $\ell_{p}\setminus P$ is discrete in $\ell_{p}$.
We start with a useful computation:

\begin{Lem} \label{liedev} Let us denote by
 $L_{ Z}$ and $L_{\overline Z}$ the Lie derivativeswith respect to the vector fields $Z$ and $\overline Z$. At all points where $Z$ is defined   
\begin{align} L_{ Z}(dd^c\r)^n=n(dd^c\r)^n=
L_{\overline Z}(dd^c\r)^n \label{lie}
\end{align}

\end{Lem}
\begin{pf} It is enough to prove the assertion on $P=\{z\in M\mid dd^c\rho(z)>0\}$,
the claim will follow by continuity. 
For any differential form $\omega$,  we recall that  
$L_Z\omega=di_Z\omega+i_Zd\omega$ where   $\iota_Z$ denotes interior multiplication. Since \((dd^c\r)^n\) is closed, \(d\r(Z)=id^{c}\r(Z)=\r\) and \(\iota_Zdd^cu=0\), using  \eqref{ddcu}, one computes
 
 \begin{align*}
 L_Z(dd^c\r)^n&=n\left(dd^c\r\right)^{n-1}\wedge L_Z(dd^c\r)\\
 &=n\left(dd^c\r\right)^{n-1}\wedge d\iota_Z(dd^c\r)\\
 &=n\left(dd^c\r\right)^{n-1}\wedge d\iota_Z(\r dd^cu + \frac{d\r\wedge d^c\r}{\r})\\
  &=n\left(dd^c\r\right)^{n-1}\wedge d\iota_Z( \frac{d\r\wedge d^c\r}{\r})\\
  &=n\left(dd^c\r\right)^{n-1}\wedge
   d\left(\iota_Z \left(\frac{d\r}{\r}\right) d^c\r- \iota_Z\left( \frac{d^c\r}{\r}\right)d\r\right)\\
 &=n(dd^c\r)^n.
 \end{align*}
The same calculation shows that 
\[
L_{\overline Z}(dd^c\r)^n=n(dd^c\r)^n
\]
 \
 \end{pf}
 
 Since \(Z(\r)=\overline{Z}(\r)\) the vector field \(\Theta = Z-\overline{Z}\) is tangent to the level set of \(\r\). Furthermore, by Lemma \ref{liedev},
\begin{align}
 L_{\Theta}(dd^c\r)^n= L_{(Z-\overline{Z)}}(dd^c\r)^n=0 \label{liedeveq}.
\end{align}

From this we get the following  key fact:

 \begin{Prop} Let  \(x\in P_{n-1}\setminus P\). Then \((dd^c\r)^n=0\) on the orbit of \(x\) in the set \(\{\r=a\}\) under the vector field  \(\Theta \). As a consequence any (extended) leaf passing through a point in $P$ lies entirely in $P$ and any leaf passing through a point in $P_{n-1}\setminus P$ lies entirely in $P_{n-1}\setminus P$.
 In particular $P$ and  $P_{n-1}\setminus P$ are foliated by leaves of the Monge-Amp\`ere foliation.
\end{Prop}
\begin{pf} 
With respect to local coordinates local coordinates $z^{1},\dots,z^{n}$ 
we have
\((dd^c\rho)^n=Ddz^1\wedge\dots d\overline z^n\), where \(D=\mbox{det} (\r_{i\overline j})\). Therefore the Lie derivative equation (\ref{liedeveq}) says that 
\[
(L_{\Theta}D)dz^1\wedge\dots d\overline z^n + DL_{\Theta}(dz^1\wedge\dots d\overline z^n)=0
\]
Since \(L_{\Theta}dz^1\wedge\dots d\overline z^n\) is an \((n,n)\) form, it is equal to \(A dz^1\wedge\dots d\overline z^n\) for some smooth function $A$ and hence \(D\) satisfies the first order linear differential equation 
\begin{align}
L_{\Theta}D+AD=0 \label{thetaeq}
\end{align}
So along an integral curve $\gamma(t)$ of \(\Theta\), equation 
\eqref{thetaeq} implies that $\tilde D(t)=_{def}D(\gamma(t))$ satisfies
\[
\frac{d\tilde D}{dt}=-A_{|\gamma(t)}\tilde D\]
and hence $\tilde D(t)\equiv 0$ if $\tilde D(0) =0$. The claim then follows.
\end{pf}

Now we are ready to prove the main result of the section:

 \begin{Prop}  Suppose that the assumptions $(A1)$, $(A2)$ and $(A3)$ hold.  Then the extended Monge-Amp\`ere foliation  $P_{n-1}$ is
holomorphic.
\end{Prop}
\begin{pf}
Let  $\ell$ be a leaf of the extended Monge-Amp\`ere foliation. At each point  $q\in \ell$ the tangent space $T_{q}\ell$ is the complex subspace of $T_{q}M$ spanned by  $Z(q)$ and hence
the leaf $\ell$ is a Riemann surface. The restriction $u_{|\ell}$  of the  function $u=\log \rho$ to $\ell$ is harmonic.  To see this, note that for $q\in M_{*}$, if $\ell$ is the leaf through $q$  and $\zeta$ is a holomorphic coordinate along $\ell$ in a coordinate neighborhood 
of $q$, then for some smooth function $\psi$ one has 
$\frac{\partial}{\partial \zeta}=\psi Z$. Hence the claim is equivalent to 
\[
dd^{c}u(\frac{\partial}{\partial \zeta},
\frac{\partial}{\partial \bar \zeta})=
\vert \psi\vert^{2}dd^{c}u(Z,\bar Z)=0
\]
where equality $dd^{c}u(Z,\bar Z)=0$  holds on $M_{*}$ as it is obvious  at  the points of the dense set $P$ where $dd^{c}\rho>0$.
\par
 Furthermore  we can see that $u_{|\ell}$ is unbounded  above on $\ell$. Suppose  in fact that $u_{|\ell}$ were bounded above and let $r =\sup_{\ell} \rho>0$. Let $p\in \overline \ell$ with $\rho(p)=r$ and let $\ell_{p}$ be the leaf of the Monge-Amp\`ere foliation passing through $p$. Then necessarily $\ell_{p}\subset \{z\in M \mid \rho(z)\geq r\}$ otherwise the leaf $\ell$ would extend past $ \{z\in M \mid \rho(z)= r\}$. But then $p$ would be a local maximum for the harmonic function $u_{|\ell_{p}}$ which is impossible since $u_{|\ell_{p}}$ is not  constant.
\par
It follows from the previous theorem that the complex gradient $Z$ extends as a real analytic vector field on
$M_{*}$ which  we denote  by $Z$. The extended vector field  $Z$  is holomorphic along any leaf of the
Monge-Amp\`ere foliation. To see this, let us cover $M_{*}$ by coordinate
neighborhoods $U$ with $C^{\infty}$-smooth coordinates
$z^{1},\dots,z^{n}$ such that the intersection of a leaf of the Monge-Amp\`ere
foliation with $U$ is given by $z^{2}=c_{2},\dots, z^{n}=c_{n}$ for suitable constants $c_{2}\dots, c_{n}$ so
that $z^{1}$ is a holomorphic coordinate along each leaf. 
Since from \eqref{Zrho} it follows that the leaves of the Monge-Amp\`ere foliation are transverse to the level set of $\rho$ and hence of $u$, on $U$ we have $du\wedge dz^{j}\neq 0$ for $j\neq 1 $ and $du=u_{1}dz^{1}+\dots+u_{n}dz^{n}$ with  $u_{1}\neq 0$. 
Furthermore since $\partial \over \partial z^{1}$ is tangent to the Monge-Amp\`ere foliation,  on  $U'=U\cap \{z\in
M_{*}\mid dd^{c}\rho(z)>0\}$ we have $Z=\varphi {\partial \over \partial
z^{1}}$ for some $C^{\infty}$-smooth function
$\varphi$. On the other hand on $U'$ we have $\rho=Z(\rho)=
\varphi {\partial \rho \over \partial z^{1}}$ from which it follows that on $U'$
we get
\begin{align}\label{extendZ two}
Z={1\over u_{1}} {\partial \over \partial z^{1}}.
\end{align}
Thus (\ref{extendZ two}) provides an expression of the $C^{\infty}$-smooth extension of $Z$ on $U$.
Since $u$ is harmonic along the leaves, we also have that   the extended vector field $Z$  is holomorphic along the intersection of $U$ with any leaf. Since $M_{*}$ is covered by such open coordinate neighborhoods, we got the claim.
\par
By continuity we notice that  $Z(\rho)=\rho$ on all $M_{*}$.
%
Then, if we denote by $X={1\over 2}(Z+ \bar Z)$  and $Y={1\over 2i}(Z- \bar Z)$  
the real and the imaginary part of  $Z$, 
we conclude that on  $M_{*}$
\begin{align}\label{XY two}
X(\rho)={1\over 4}\rho\hskip1cm{\rm and}\hskip1cm Y(\rho)=0.
\end{align}
Since $\rho\colon M_{*}\to (0,+\infty)$ is proper and the level sets of $\rho$
are compact, it follows from (\ref{XY two}) that $X$ and $Y$ are complete i.e.
the flows $\phi$ and $\psi$ of $X$ and $Y$ respectively, are both defined on
$\R$.
If $l_{p}$ is the leaf through any fixed point  $p\in M_{*}$, then the map
$f \colon \C \to l_{p}$ defined by
\[f(t+is)=\phi(t,(\psi(s,p))\]
is holomorphic since $f'(t+is)=Z(f(t+is))$, $Z$ is holomorphic along the leaf and $f$ is non degenerated as $Z\neq0$ on $M_{*}$. Therefore the leaf $l_{p}$ is a
parabolic Riemann surface. Since this is true for all leaves, Burns's Theorem  $3.2$ in
\cite{Burns} applies on the set $P$ and implies that the vector field $Z$ is holomorphic on $P$. As $Z$ is smooth on $P_{n-1}$, 
the extended Monge-Amp\`ere  foliation is holomorphic.
\end{pf}

\section {Extending the foliation on higher degeneracy set}
\setcounter{equation}{0}

In this section we shall prove that indeed the Monge-Amp\`ere foliation extends  to the set where the rank of the form $dd^{c}\rho$ is less than $n-1$. To this purpose we use a remark that, again, we can  prove  only assuming real analyticity:

\begin {Lem} \label{weak} Suppose that the assumptions $(A1)$, $(A2)$ and $(A3)$ hold. The set of ``weak points''
$$W=\{p\in M_{*} \mid {\rm rank }dd^{c}\rho \leq n-2\}$$
 does not contain   divisors (i.e.   complex hypersurfaces)
\end{Lem}
\begin{pf} Suppose that $S=\{\rho=c>0\}$ is a generic level set of 
$\rho$ which intersects $W$. Then $S\cap W$ is a real analytic subvariety of $S$. For a (local) real analytic  subvariety $V$ of $S$, 
let us denote ${ hol dim} (V)$ the {\sl holomorphic dimension} of $V$ as defined in \cite{DiederichFornaess}, Definition 1, pg 372 (see also
\cite{Kohn2} where the notion was originally defined).
We claim that
$${ hol dim} (S\cap W) \geq 1.$$
This fact is consequence of the following observation: at a point $z\in S\cap W$ the Levi form of $\rho$ has 
 eigenvalue $0$  with multiplicity at least $2$ and hence there is
 at least $1$ eigenvector with eigenvalue  $0$ which is in the intersection of the holomorphic tangent space to $S$ and the complex tangent space to W at $z$.
 Then by Theorem 3 of  \cite{DiederichFornaess} there is a 
 complex curve (i.e. complex submanifold of dimension $\geq 1$) contained in $S$. But this cannot happen since $S$ is a 
 compact real analytic hypersurface (see Theorem 4 of  \cite{DiederichFornaess}). 
 \end{pf}

Using Lemma \ref{weak} we have the desired extension result
for the vector field $Z$ and the Monge-Amp\`ere foliation:
\noindent

\begin {Prop}  The vector field
$Z$ and the Monge-Amp\`ere foliation extend holomophically to the set of ``weak points'' $W$.
\end{Prop}
\begin{pf} The vector field
$Z$ is holomorphic on $M_{*}\setminus W$. Since $W$ does not contain a complex hypersurface, it follows that $W$ is a removable 
set for the vector field $Z$ (see  \cite{merker}, Theorem 2.30, Chapter VI, and its proof). The extended vector field, which we still denote 
$Z$, is holomorphic and nonzero on $M_{*}$ since by continuity
$Z(\rho)=\rho$. Furthermore, again by continuity, we have
$dd^{c}\rho(Z, \overline Z)=0$ also on $W$. Therefore we get the 
extension and  the holomorphicity of the Monge-Amp\`ere foliation
on all $M_{*}$.
 \end{pf}
 
Our main result can now be stated:

\begin{Th}  \label{structure}
 Let  $M$ be a Stein manifold of dimension $n$ and 
 $\rho \colon M \to [0,+\infty)$  be a
continuous plurisubharmonic exhaustion of class $C^\omega$ with   $d\rho \not= 0$ on $M_{*}=\{\rho >0\}$
and such that
\\
$(i)$  the open set  $M_{*}$  is connected,
\\ 
$(ii)$ on \(M_{*}\)  the function $u=\log \rho$ satisfies the complex
homogeneous Monge-Amp\`ere equation
\begin{align}(dd^c u)^n = 0,
\end{align}
Then the  complex gradient vector field $Z$  of $\rho$ defined 
by \eqref{zeta}
on $P=\{p\in M_{*}\mid (dd^{c}\rho)^{n}\neq0\}$,
extends as a  holomorphic vector
field on $M$ and the minimal set $M_{0}=\{z\in M \mid \rho(z)=0\}$ reduces to one  point.
\end{Th} 
\begin{pf} 
We only need to show that the minimal set $V$ of
$\rho$ reduces to one point. We observe that \(M_{0}\) cannot be empty as
otherwise $u=\log\rho$ would be a smooth exhaustion which solves the complex
homogeneous Monge-Amp\`ere equation on a Stein manifold. But such  function do
not exist (Theorem 1.1 of \cite{lempert-szoke}). The holomorphic flow
$\Psi\colon \C\times M_{*}\to M_{*}$ on $M_{*}$ of the vector field $Z$ extends
to a holomorphic flow $\Psi\colon \C\times M\to M$ as $M_{0}$ is compact. On the other hand for any $p\in M_{0}$, since the flow of any point in $M_{*}$ is contained in $M_{*}$ we
must have
$\Psi(\C\times\{p\})\subset M_{0}$, which, by Liouville
theorem
implies that $\Psi(\C\times\{p\})=\{p\}$. Thus $Z(p)=0$. But then
$M_{0}=\{z\in M \mid Z(z)=0\}$ is a compact analytic set,
i.e. a finite set of points. To finish our proof it is enough to show that the
minimal set of $\rho$ is connected. This can be done by repeating verbatim the
argument
at page 357 of \cite{Burns}. We give here an outline of the proof.  Suppose that
for some  non empty disjoint compact subsets $K_{1}, K_{2}\subset M$ we
have $M_{0}=K_{1}\cup K_{2}$. Let $V_{1}, V_{2}$ disjoint
open sets with $K_{i}\subset V_{i}$ for $i=1,2$.  For $r>0$ sufficiently small
we have $S_{i}=\rho^{-1}(r)\cap V_{i}\neq \emptyset$ for $i=1,2$.
If $G\colon \R\times M_{*}\to  M_{*}$ is the flow of the real part $X$ of the
vector field $Z$, then $G(\R\times\rho^{-1}(r))=M_{*}$. On the other hand
if $U_{i}=G(\R\times S_{i})$ for $i=1,2$, then $U_{1},U_{2}$ are open disjoint
subsets with $U_{1}\cup U_{2}=M_{*}$ contradicting the fact that $M_{*}$ is
connected.

\end{pf}

\section {Applications and final remarks}
\setcounter{equation}{0}

As announced, Theorem \ref{structure} can be used to give a full answer to a question of D. Burns (\cite{Burns}):

\begin{Th} Let $\r$ be a  positive homogeneous polynomial on
${\bf C}^n$  such that $u=\log \r$ is plurisubharmonic and satisfies
\begin{align} \label{MA}(dd^c u)^n = 0.
\end{align}
Then  $\r$ is a homogeneous polynomial of bidegree $(k,k)$.
\end{Th}

\begin{pf} 
Using Theorem \ref{structure}, the proof is a repetition of the elementary argument used to prove   Theorem 3.1 of \cite{Kalka-Patrizio2}. We give here a very short outline of the key steps of the proof for completenss referring the reader to \cite{Kalka-Patrizio2} for the details.   Theorem \ref{structure} implies that the complex gradient $Z$, defined  on $P$ by
\begin{align}\label{zedbis}
Z=\sum_{\mu}Z^{\mu}{\frac{\partial\,} {\partial z^{\mu}}} =
\sum_{\mu,\nu}Z^{\mu}\rho^{\mu\bar{\nu}}\rho_{\bar\nu}\frac{\partial}{\partial
z^{\mu}},
 \end{align}
extends holomorphically to all ${\bf C}^n$. On the other hand,  (\ref{zedbis})
shows that $Z$ is homogeneous of degree one on a dense subset of
${\bf C}^n$, so that $Z$ is   in fact linear.
If
$\displaystyle{\r=\sum_{l+m=2k}\r^{l,m}}$
is the decomposition of $\r$ in homogeneous polynomial of bidegree $(l,m)$, a bidegree
argument using the fact that on $P$ one has $Z(\rho)=\rho$, shows that 
$$0=\r^{0,2k}_{\bar\alpha}=\r^{2k,0}_{\bar\alpha} \hskip1cm
{\rm and}\hskip1cm
\r^{l,m}_{\bar\alpha}=\sum_{\mu}Z^{\mu}\r^{l,m}_{\mu\bar\alpha}
$$
for every $l,m$ with $l+m=2k, l,m\geq 1$.
If $w=(w^1,\dots,w^n)\in P$, then $dd^c\r^{k,k}_{|w}>0$
 and, again by bidegree reasons and the fact that $\r^{k,k}$ is homogeneous of bidegree $(k,k)$,
then
\begin{align}\label{newzed} 
Z^{\mu}_{|w}=
(\r^{k,k})^{\bar\alpha\mu}(w)\r^{k,k}_{\bar\alpha}(w)={1\over k}w^\mu.
\end{align}
By continuity (\ref{newzed}) holds on all ${\bf C}^n$.
As consequence the (extended) Monge-Amp\`ere foliation associated to $u=\log\r$ is given by the foliation of  ${\bf C}^n\setminus\{0\}$ in lines through the origin and the restriction of $\log \r$ to any complex
line through the origin is harmonic with a logarithmic singularity of weight $2k$ at the origin.
Thus, for any $0\not=z\in{\bf C}^n$  and $\lambda \in {\bf C}$,
$$\log \r(\lambda z)=2k\log \vert \lambda \vert + O(1).$$
Hence the restriction of $\r$ to any complex line through the origin is
homogeneous of bidegree $(k,k)$ and therefore $\r$ is a homogeneous polynomial
of bidegree $(k,k)$ on
${\bf C}^n$.

\end{pf}

Again using Theorem \ref{structure} and, with the obvious formal changes,  the arguments outlined in section 4 of \cite{Kalka-Patrizio2}, we have the following classification result:

\begin{Th} Let $M$ be a Stein manifold of dimension $n$ equipped 
with a continuous  plurisubharmonic exhaustion 
 $\rho \colon M \to [0,+\infty)$
of class $C^\omega$ with   $d\rho \not= 0$ on $M_{*}=\{\rho >0\}$
and such that
\\
$(i)$  the open set  $M_{*}=\{\rho >0\}$  is connected,
\\
$(ii)$  on $\{\rho >0\}$  the function $u=\log \rho$ satisfies the complex
homogeneous Monge-Amp\`ere equation,
\begin{align}(dd^c u)^n = 0.
\end{align}
Then there exists a biholomorphic map $\Phi\colon {\bf C}^n\to M$,
such that, for suitable positive real numbers $c_{1},\dots,c_{n}$, the pull back
$\r_{0}=\r\circ\Phi$ of the exhaustion $\r$ satisfies
\begin{align}\label{rholinerized}
\rho_{0}(e^{c_{1}\lambda}z^{1},\dots,e^{c_{n}\lambda}z^{n})=
\vert e^{\lambda}\vert^{2}\rho_{0}(z).
\end{align}
so that the sublevelsets of $\rho$ are biholomorphic to generalized weighted
circular domains.
\end{Th}

\begin{pf} Also in this case  Theorem \ref{structure} provides the necessary tools to adapt at this situation the arguments of (\cite{Wong}). In fact Theorem \ref{structure} shows that  the
Monge-Amp\`ere foliation associated to $u=\log \rho$ extends to a holomorphic foliation throughout $M_{*}=\{\rho >0\}$ and the complex gradient vector field $Z$, locally defined by
\begin{align}\label{zed}
Z=\rho^{\mu\bar{\nu}}\rho_{\bar\nu}\frac{\partial}{\partial
z^{\mu}},
 \end{align}
on the set  $P=\{dd^{c}\rho>0\}$, extends holomorphically to all $M$ so that the extension, which we shall keep denoting $Z$, is tangent to the leaves of the extended Monge-Amp\`ere foliation  on $M_{*}$ and
the equality $Z(\rho)=\rho$,
which is equivalent to the Monge-Amp\`ere equation  on the set
where $P$, holds, by continuity on all $M$.
Finally we have that the minimal set of the function $\rho$ reduces to a single point:
$\{\rho =0\}=\{O\}$. We shall call the point $O$ the {\sl center} of $M$.

The rest of the proof  in  (\cite{Wong}) uses only these facts and no other consequence of strict pseudoconvexity. Exactly as for the proof Theorem 3.1 in \cite{Kalka-Patrizio2}, it is enough
to outline the main steps of the rest of the proof there to
underline the minor variations needed under our weaker assumptions.
\\ \\
{\bf Step 1:} For any $r_{1}, r_{2}>0$ the level sets $\{\rho=r_{1}\}$
and $\{\rho=r_{2}\}$ are $CR$ isomorphic and the sublevel sets $\{\rho<r_{1}\}$
and $\{\rho<r_{2}\}$ are biholomorphic.
\\ \\
This follows from a Morse Theory type of argument. The flow of the real part $X=Z+\bar Z$ of the (extended) complex gradient vector field maps level defines a local group of biholomorphims when $Z$ is holomorphic and maps level sets of $\r$ into other level sets. Details of the argument can be found on page 24 of (\cite{Patrizio-Wong}). There it was assumed   $\r$ to be strictly plurisubharmonic  merely to ensure that the vector field $Z$ is defined everywhere which in our case we prove by other means. The claim is an application of Bochner-Hartogs extension theorem.
\\ \\
{\bf Step 2:} There exists a coordinate neighborhood $U$ of the center $O$, with coordinates $z^{1},\dots z^{n}$ centered at  $O$, and positive real numbers $c_{1},\dots, c_{n}>0$, on $U$, the vector field $Z$ has the following
expression:
\begin{align}\label{zedlinerized}
Z=c_{1}z^{1}{\partial \over \partial z^{1}}+\dots
c_{n}z^{n}{\partial \over \partial z^{n}}. \end{align}
\\
The expression for the vector field $Z$ can be proved as on
on page 248 of \cite{Wong} as a consequence of a standard linearization theorem (see  \cite{Bochner-Martin}, Chapter III, for instance)
since the equation $Z(\rho)=\rho $ implies that the one parameter group
associated with imaginary part of $Z$ preserves the level sets of $\rho$ and
therefore has compact closure in the automorphism group of each sublevel set.
\\ \\
{\bf Step 3:} There exist  suitable positive real numbers $c_{1},c_{2}$ so that
for $z=(z^{1},\dots, z^{n})\in U$ and $\lambda \in \C$ so that
$(e^{c_{1}\lambda}z^{1},\dots, e^{c_{n}\lambda}z^{n})\in U$, we have:
\begin{align}\label{rholinerized}
\rho(e^{c_{1}\lambda}z^{1},\dots, e^{c_{n}\lambda}z^{n})=
\vert e^{\lambda}\vert^{2}\rho(z).
\end{align}
\\
As in page 249 of \cite{Wong} one integrates
vector field $Z$ explicitly. The only relevant fact needed here is that
$Z(\rho)=\rho$.
\\ \\
What is left to prove is that there is a global biholomorphism of
${\bf C}^n$ onto $M$.   Step 3 implies that for  $\epsilon>0$ small enough there exists a biholomorphism
$\varphi : M(\epsilon)\to G(\epsilon)$ of the
sublevel set $M(\epsilon)$ of $\rho$  into a (fixed) weighted
circular domain
$G(\epsilon)=\{z\in  {\bf C}^n \mid \r_{0}(z)<\epsilon\}$
Then
$\r_{0}=\r\circ \varphi$ on $G(\epsilon)$  is defined on all ${\bf C}^n$ and
it satisfies $(i)$, $(ii)$, $(iii)$.
By Step 1 it follows that the flows of the real parts of the complex gradients
of $\r$ and $\r_{0}$ are biholomorhisms of sublevel sets.
The required biholomorphic map $\Phi \colon {\bf C}^n \to M$ can
be defined by composition of $\varphi$ and the flows of the real parts
of the complex gradients.

\end{pf} 

A final remark is in order. In complex dimension $n=2$ we were able to obtain the main results of this paper under the hypothesis that the 
exhaustion is smooth and satisfies a finite type condition. This was possible because of  results which hold only for dimension $n=2$ -- in particular the fact that foliation with parabolic leaves are always holomorphic in that case -- and the fact that, simply for dimensional reasons, the Levi form of the solution of the Monge-Amp\`ere restricted to the holomorphic tangent space to a level set may have at most one zero eigendirection. 
In higher dimension, of course, higher degeneracies are possible and we are able to overcome this difficulty only using in an essential 
way results which need real analyticity due to Burns and Bedford 
(\cite{Bedford-Burns}) and Diederich and Fornaess 
(\cite{DiederichFornaess}). Whether or not it is possible to replace
real analyticity with assumptions such as finite type, or the like,
remains open.

\bigskip
\bigskip
\font\smallsmc = cmcsc8
\font\smalltt = cmtt8
\font\smallit = cmti8
\hbox{\parindent=0pt\parskip=0pt
\vbox{\baselineskip 9.5 pt \hsize=3.1truein
\obeylines
{\smallsmc
Giorgio Patrizio
Dip. Matematica ``U. Dini''
Universit\`a di Firenze
Viale Morgani 67/a
I-50134 Firenze
ITALY
}\medskip
{\smallit E-mail}\/: {\smalltt patrizio@math.unifi.it
}
}
\hskip 0.0truecm
\vbox{\baselineskip 9.5 pt \hsize=3.7truein
\obeylines
{\smallsmc
Morris Kalka
Department of Mathematics
Tulane University
New Orleans, LA  70118
USA
}\medskip
{\smallit E-mail}\/: {\smalltt kalka@math.tulane.edu
}
}
}

\end{document}